\documentclass{amsart}
\usepackage{epsfig}

\theoremstyle{plain}
\newtheorem{theorem}{Theorem}[section]

\theoremstyle{definition}

\theoremstyle{remark}
\newtheorem{remark}[theorem]{Remark}

\begin{document}

\title[Twisted torus knots $T(p,q,3,s)$ are tunnel number one]
{Twisted torus knots $T(p,q,3,s)$ are tunnel number one}

\author{Jung Hoon Lee}
\address{School of Mathematics, KIAS\\
207-43, Cheongnyangni 2-dong, Dongdaemun-gu\\
Seoul, Korea. Tel:\,+82-2-958-3736}
 \email{jhlee@kias.re.kr}
%\thanks{}


\begin{abstract}
We show that twisted torus knots $T(p,q,3,s)$ are tunnel number one.
A short spanning arc connecting two adjacent twisted strands
is an unknotting tunnel.
\end{abstract}

\maketitle

\section{Introduction}

A {\it torus knot} is a knot that can be embedded in a standard torus in $S^3$.
For a fixed standard meridian and longitude system on a standard torus,
let $T(p,q)$ denote a torus knot that runs $p$ times in longitudinal direction and $q$ times
in meridional direction. By convention, we assume that $p>q>0$ and $p$ and $q$ are relatively prime. (If $p<0$ or $q<0$, we consider it as a mirror image.)

Take $r$ ($1<r<p$) adjacent parallel strands of $T(p,q)$ and replace them with $s$ times full twists. The resulting knot is called a {\it twisted torus knot} $T(p,q,r,s)$.
The class of twisted torus knots is interesting in many aspects. For example, it gives strong candidates for non-minimal genus, weakly reducible and unstabilized Heegaard splittings obtained by boundary stabilization \cite{Moriah-Sedgwick}.

A knot $K$ in $S^3$ is said to be {\it tunnel number one} if there exists an arc $\gamma$ properly embedded in the exterior of $K$ such that $cl(S^3-N(K\cup\gamma))$ is a genus two handlebody. It is well known that torus knots and $T(p,q,2,s)$ are tunnel number one. It is also known that tunnel number of a twisted torus knot is one or two. As a next step to $T(p,q,2,s)$, Morimoto asked about the tunnel number of $T(p,q,3,s)$ (\cite{Morimoto}, Problem $5'$). We answer the question.

\begin{theorem}
Twisted torus knots $T(p,q,3,s)$ are tunnel number one.
\end{theorem}

As a corollary, every $T(p,q,3,s)$ is a prime knot since tunnel number one knots are prime \cite{Norwood}, \cite{Scharlemann}.

\begin{remark}
There exist composite twisted torus knots \cite{Morimoto}.
\end{remark}

\section{Proof of Theorem 1.1}

\begin{figure}[h]
   \centerline{\includegraphics[width=12cm]{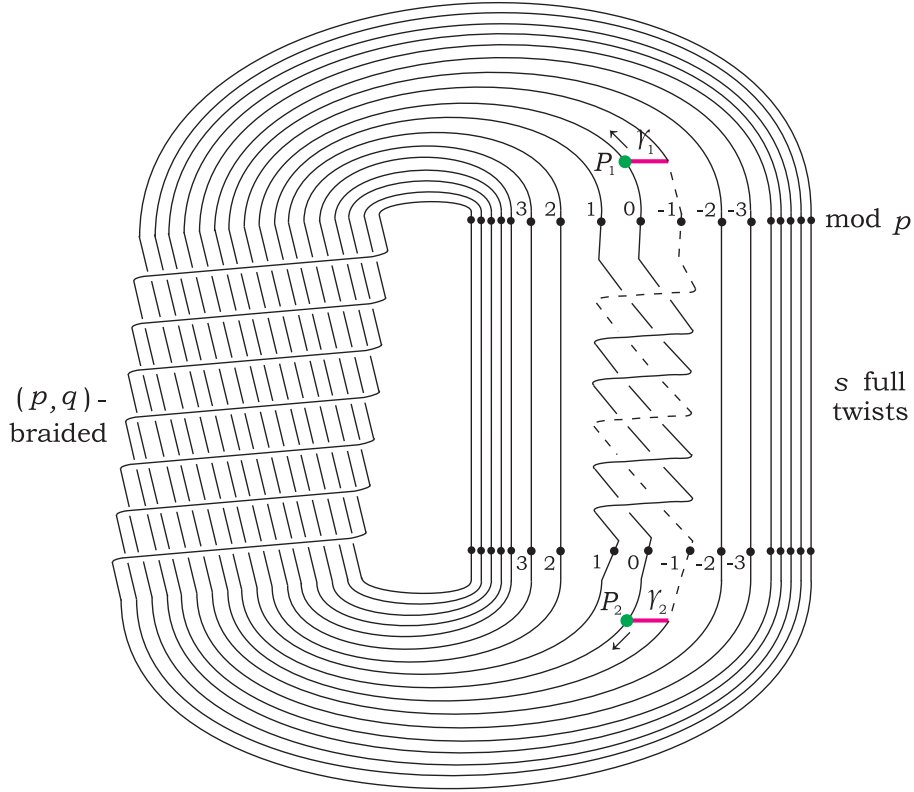}}
    \caption{$T(17,7,3,-2)$}
\end{figure}

Consider a standard knot diagram of a twisted torus knot $T(p,q,3,s)$ as in the Figure 1. (It is $T(17,7,3,-2)$.)
There is $(p,q)$-braided part on the left, and $s$ times full-twisted part on the right. On the full-twisted part of the diagram, label the middle strand among the three full-twisted strands as $0$. Label all the other strands with consecutive integers modulo $p$.

Connect two strands labelled as $0$ and $-1$ with a short spanning arc $\gamma$. It will be clear later that $\gamma$ is an unknotting tunnel. Do slide and isotopy through the full-twisted part so that the twisted part becomes $s$ times full twists on two strands. (The dotted arc in Figure 1. will disappear after the slide and isotopy.) Now we have two copies of $\gamma$, $\gamma_1$ and $\gamma_2$, where $\gamma_1$ is above
the twisted part and $\gamma_2$ is below the twisted part.

Let $P_1$ and $P_2$ be endpoints of $\gamma_1$ and $\gamma_2$ on the strand labelled $0$, respectively.
Suppose $P_1$ moves along $T(p,q,3,s)$. Starting to the upper direction, suppose $P_1$ arrives at the strand labelled $1$ before arriving at the strands labelled $2$ or $-1$ in the bottom part of the twists. That implies $\gamma_1$ is isotopic to an arc connecting two parallel strands of $s$ times full twists, which is an unknotting tunnel for the type $T(-,-,2,s)$ and the proof is done. The same arguments holds for $P_2$. So we only need to show that either $P_1$ or $P_2$ arrives at the strand labelled $1$ before arriving at the strands labelled $2$ or $-1$.

Suppose it does not happen. Then one of the following properties holds.

\begin{itemize}
\item[(*)] $P_1$ arrives at the strand $2$ first and $P_2$ arrives at the strand $-1$ first among $\{1,2,-1\}$, or

\item[(**)] $P_1$ arrives at the strand $-1$ first and $P_2$ arrives at the strand $2$ first among $\{1,2,-1\}$
\end{itemize}
\vspace{0.1cm}

Case $1$)\, Suppose (*) holds. When $P_1$ passes through the $(p,q)$-braided part from above to below, the labelled number on the strand it belongs increases by $q$ modulo $p$. For $P_2$, it decreases by $q$ modulo $p$. Hence by property $(*)$, there exist some $k$ ($k<p$) and $j$ ($j<p$) satisfying the following equations.
\begin{equation}
q\neq 1,-1,\quad 2q\neq 1,-1\quad \ldots\quad kq=2 \pmod{p}
\end{equation}
\begin{equation}
-q\neq 1,2,\quad -2q\neq 1,2\quad \ldots\quad -jq=-1 \pmod{p}
\end{equation}

So $jq=1 \pmod{p}$ and we can see that $j>k$ from (1). Adding equations, we get $-(j-k)q=1 \pmod{p}$ and this contradicts equations (2).

\vspace{0.1cm}
Case $2$)\, Suppose (**) holds. There exist some $k$ ($k<p$) and $j$ ($j<p$) satisfying the following equations.

\begin{equation}
q\neq 1,2,\quad 2q\neq 1,2\quad \ldots\quad kq=-1 \pmod{p}
\end{equation}
\begin{equation}
-q\neq 1,-1,\quad -2q\neq 1,-1\quad \ldots\quad -jq=2 \pmod{p}
\end{equation}

So $-kq=1 \pmod{p}$ and we can see that $k>j$ from (4). Adding equations, we get $(k-j)q=1 \pmod{p}$ and this contradicts equations (3).

\vspace{0.1cm}
 The sign of $s$ does not effect much and the arguments are similar. This completes the proof.

\end{document}